\input amstex.tex
\input amsppt.sty

\def\Sup{{\roman{Supp}}}

\magnification=1200
\tolerance=5000
\refstyle{A}

\topmatter
\title Split Clifford modules over a Hilbert space
\endtitle
\thanks
This research was supported in part by CONICET, FONCYT, CONICOR and UNC.
\endthanks
\author Esther Galina, Aroldo Kaplan and Linda Saal
\endauthor
\affil
CIEM-FaMAF, Universidad Nacional de C\'ordoba\\
\endaffil
\email
galina\@mate.uncor.edu \\
kaplan\@mate.uncor.edu \\
saal\@mate.uncor.edu
\endemail
\address
CIEM - FAMAF, Universidad Nacional de C\'ordoba, Ciudad Universitaria,
(5000) C\'ordoba, Argentina
\endaddress

\

\abstract {We describe the real, complete modules of the Clifford Algebra $C(Z)$, where $Z$ is a real Hilbert space acting by skew-adjoint operators.}
\endabstract
\endtopmatter

\document

\subhead 1. Introduction
\endsubhead

 Let $Z$ be a real separable prehilbert space and  $C(Z)$ the corresponding Clifford Algebra, i.e.,  the quotient of the Tensor Algebra of $Z$ by the ideal generated by the elements of the form
 $z_1z_2+z_2z_1 + 2(z_1,z_2)1$, with $z_1,z_2\in Z$. 
We determine here all the representations of $C(Z)$ on real separable Hilbert spaces where the elements of
$Z\subset C(Z)$ act skew-adjointly. Equivalently, all sequences $J_1,J_2,...$ of orthogonal, mutually anticommuting complex structures on a real Hilbert space. 

Without loss of generality we may assume that $Z$ is complete as well.
Every such representation gives rise, upon complexification, to a representation of the algebra $C_\Bbb C(Z)={\Bbb C}\otimes C(Z)$ on a complex Hilbert space, whith $Z$  acting by skew-hermitian operators. 
A complete set of representatives of such representations up to equivalence was given by 
G\aa rding and Wightman \cite {GW}. 

In this article we determine which G\aa rding-Wightman modules split over $\Bbb R$ and describe all their $C(Z)$-invariant real forms up to equivalence. These effectively solves the problem, since any $C(Z)$-module must arise as one of these invariant forms.

To put the results in perspective, recall that
for $m=\dim Z <\infty$, the question of when a given $C_\Bbb C(\Bbb R^m)$-module, or space of spinors, splits over $\Bbb R$, goes back to E. Cartan and Killing [C], who also provided the answer and geometric interpretation. The representations of $C_\Bbb C(\Bbb R^m)$ are direct sums of either one or two finite dimensional irreducible representations, depending on whether $m$ is even or odd. These remain irreducible over $\Bbb R$ if $m\equiv 1,2,$ mod. (4) and split into a direct sum two dual irreducibles otherwise. In particular, {\it a $C_\Bbb C(\Bbb R^m)$-module splits over $\Bbb R$ if and only if} $m\equiv 0, 3$ mod.(4). We recover this result, at least when $m$ is even.

If $\dim Z=\infty$ the situation is much more complicated, as can be expected. Rather, the remarkable feature of the infinite Clifford Algebra -as well as its bosonic cousin, the Heisenberg Lie Algebra, is that its representations can be parametrized at all and the possible real forms described. In \cite {GW}, the  resulting set is described as ``a true maze". Since both there and here the matters of equivalence and irreducibility are not completely settled, an alternative argument is required to analize the splitting problem in infinite dimensions.
In any case, although  Clifford and Heisenberg Algebras are exceptions to Kirillov's dictum that ``in infinite dimensions there are no theorems, only examples", it is the examples what may matter most.  

The basic Fermi-Fock representations, familiar to physicists, do not split and in fact are irreducible over $\Bbb R$. But we find others that are irreducible over $\Bbb C$, split over $\Bbb R$ and the corresponding factors are essentially new, i.e., they do not arise by restriction of any $C_\Bbb C(Z)$-module.

Although the examples of the latter kind may fail as models for quantum fields because they lack a vacuum (or only apply to systems with infinite energy), they do have natural realizations in $L_2({\Bbb T})$ which are closely related to wavelets and multilinear operators. The Clifford Algebra acts by singular integral transforms on the circle and the corresponding splittings can be interpreted in terms of wavelets. Conversely, this picture describes wavelets as infinite-dimensional spinors.

We are grateful to H. Araki, A. Grunbaum, A. Jaffee, V. Kac and J. Vargas for their useful advise.

\subhead 2. The G\aa rding-Wightman modules\endsubhead

We describe here the results of [GW] and translate them in the language of Clifford modules (see also [G], [BSZ]).
According to [GW], a representation of the Canonical Anticommutation Relations (CAR) consists of a countable set of bounded linear operators $a_k$ acting on a separable complex Hilbert space $H$, satisfying 
$$
a_ja_k + a_ka_j =0 \qquad \qquad a_ja^*_k + a^*_ka_j = \delta_{jk}.
$$
The $a_k$ and $a_k^*$ are sometimes called the (fermionic) creation and anihilation operators, respectively.

To view these as representations of $C_\Bbb C(Z)$, we choose a complex structure $\sigma$ on $Z$ compatible with its metric, making $(Z,\sigma)$ into a complex-hermitian Hilbert space and a unitary basis $z_1, z_2, ...$ of it, so that $\{z_k, \sigma(z_k)\},$ is a real orthonormal basis of $Z$. Then, the following formulas define a 1-1 correspondence between representations of the CAR and continuous representations $\Bbb J$ of the algebra $C_\Bbb C(Z)$ on $H$ such that the operators $\Bbb J_z$, $z\in Z$, are skew-hermitian:
$$
a_k=\frac 1{ 2} \left(\Bbb J_{\sigma(z_k)} + i \Bbb J_{z_k}\right) \qquad \qquad
a^*_k=\frac 1{ 2} \left(-\Bbb J_{\sigma(z_k)} + i \Bbb J_{z_k}\right).
\tag {2.1}
$$
Note that $z\mapsto \Bbb J_z$ is not assumed to be $\Bbb C$-linear and that for $z\in Z$, $|z|=1$, the operator $\Bbb J_z$ is unitary and its square is $-I$.

Let $X$ be the set of sequences $x=(x_1,x_2,\dots)$ where each $x_i$ is $0$ or $1$ (space of fermionic ``occupation numbers"). $X$ is a group under componentwise addition modulo $2$. Let $\Delta$ be the subgroup of $X$ consisting of sequences with only finitely many $1$'s. It is generated by the sequences $\delta_k$ where the $k$-th component is $1$ and the all others are $0$. 
The G\aa rding-Wightman (or GW) modules are parametrized by triples 
$$(\mu,\nu,\{c_k\})\tag {2.2}$$ with the following objects as entries:
\smallskip

$\bullet\ \ \mu$ is a positive  measure over the $\sigma$-algebra generated by the sets $X_k=\{x\: x_k=1\}$, which is quasi-invariant under translations by $\Delta$. This means that its translates
 $\mu_\delta(E)=\mu(E+\delta),$ $\delta \in \Delta$,  all have the same sets of measure zero and implies the existence of the Radon-Nikodym derivatives $d\mu_\delta/d\mu$.
\smallskip

$\bullet\ \ \nu: X\rightarrow \Bbb Z_+ \cup \{\infty\}$ is a measurable function such that $\nu(x+\delta) =\nu(x)$ for all $\delta\in \Delta$ and almost all $x\in X$. 
\smallskip

$\bullet\ \ c_k(x): H_x \rightarrow H_{x+\delta_k}\approx H_x$ is a sequence of unitary operators depending measurably on $x$ and satisfying 
$$
\aligned
c^*_k(x)&=c_k(x+\delta_k)\\
c_k(x)c_l(x+\delta_k)&= c_l(x)c_k(x+\delta_l)
\endaligned
\tag {2.3}
$$
The Hilbert space isomorphisms $H_{x+\delta_k}\approx H_x$ are part of the data, although changing them results in unitarely equivalent representations. One may therefore assume that the family $H_{x}$ is invariant under translations by $\Delta$.

The representation associated to a triple 
$(\mu,\nu,\{c_k\})$ is realized in the Hilbert space 
$$H=\int_X^\oplus H_x \ d\mu(x)
\tag {2.4}
$$
and defined by the following formulae, where an $f\in H$ is regarded as an assignement $x\mapsto f(x)\in H_x$:
$$
\aligned
(\Bbb J_{z_k} f)(x)&= (-1)^{x_1+\dots +x_{k-1}+1}i c_k(x) \sqrt{\frac {d\mu(x+\delta_k)}{d\mu(x)}} f(x+\delta_k)\\
(\Bbb J_{\sigma(z_k)} f)(x)&= (-1)^{x_k}i(\Bbb J_{z_k} f)(x)\\
\endaligned
\tag {2.5}
$$

The main observation here is that in any representation of the CAR, the operators $N_k= a^*_k a_k$ and $N'_k=a_ka^*_k$ form a conmuting set of projections over the sets $X_k$ and their complements $X'_k$, respectively, so that
$$
(N_k f)(x)=x_k f(x) \qquad \qquad (N'_k f)(x)=(1-x_k) f(x)
\tag {2.6}
$$
for all $f\in H$ and for almost all $x$. The direct integral (2.4) is just the corresponding spectral decomposition. 

Setting
 $$A_k = i\Bbb J_{z_k}\qquad B_k =  -i\Bbb J_{\sigma(z_k)}, $$
our notation becomes consistent with that of [GW]. 

\

\subhead 3. Real Clifford modules
\endsubhead

If $V$ is a (real) module over $C(Z)$, then $\Bbb C\otimes V$ is a (complex) module over $\Bbb C\otimes C(Z)$ with the 
$C(Z)$-invariant decomposition over $\Bbb R$
$$\Bbb C\otimes V = V \oplus iV;$$
$V$ is an {\it invariant real form} of $\Bbb C\otimes V$. Hence, to determine all real modules over $C(Z)$ it is enough to determine all the invariant real forms of the G\aa rding-Wightman modules. This is equivalent to determining all the antilinear operators $S:H\to H$ (the corresponding complex-conjugations) which commute with $C(Z)$ and such that
$$
S^2=1, \qquad ||Sf|| = ||f||. \tag {3.1}
$$
As $S$ is  a $C(Z)$-morphism and $A_k = i \Bbb J_{z_k}$ and $B_k = - i \Bbb J_{\sigma(z_k)}$ by (2.1) and (2.5), we know that for all $k$
$$
S A_k= -A_k S  \qquad \text{and} \qquad  S B_k=- B_k S.
\tag {3.2}
$$
Hence,
$$
S N_k=N'_k S.
\tag {3.3}
$$
Every operator $L$ that conmutes with $N_k$ and $N'_k$ for all $k$ is an  operator of multiplication by an essentially bounded measurable operator-valued function $\phi$: 
$$(Lf)(x)=(L_\phi f)(x)= \phi(x)f(x)$$
(cf. \cite G ). In particular, $N_k$ and $N'_k$ correspond to multiplication by the characteristic functions of the sets $X_k$ and $X'_k$ respectively (see (2.6)).

From now on, we let $H$ be the complex Clifford module associated to $(\mu,\nu,c_k)$ as in the previous section.

\proclaim {Lemma 3.4} Suppose that $H$ admits an invariant real structure $S$. Let $E\subset X$ be measurable and  $\chi_E$ its characteristic function. Then, if 
$$\bold 1 - E = \{\bold 1-x=(1-x_1,1-x_2,...)\: x\in E\},
$$

(i) $SL_{\chi_E} =L_{\chi_{\bold 1 - E}} S$.

(ii) $\Sup (Sf) = \bold 1-\Sup (f)$.

(iii) $(Sf,Sg) = \overline{(f,g)}$.
\endproclaim

\demo{Proof} $(i)$ follows from (3.3) and the fact that $X_k$ and $X'_k $ generate the $\sigma$-algebra where $\mu$ is defined.
To prove $(ii)$, set $F=\Sup(f)$. Then
$$
\Sup(Sf)=\Sup(S\chi_F f) = \Sup(\chi_{\bold1-F}Sf) \subset \bold1-F.
$$
As $S$ is an involution,
$$
F=\Sup(S^2 f)  \subset \bold1- \Sup(Sf)\subset \bold1 - (\bold1-F)=F.
$$
Consequently, all these inclusions are equalities and the assertion follows.
For $(iii)$, as $S$ preserves norm, ${\roman Re}(Sf,Sg)={\roman Re}(f,g)$. Also, 
$${\roman Im}(Sf,Sg) = -{\roman Re}(iSf,Sg) = 
{\roman Re}(Sif,Sg)= {\roman Re}(if,g) = 
-{\roman Im}(f,g),$$
finishing the proof.
\qed
\enddemo

If $\mu$ is a measure on $X$, then
$$
\tilde \mu(E) = \mu( \bold 1 - E)
$$
is another.

\proclaim{Theorem 3.5} $H$ admits an invariant real form if and only if $\mu$ and $\tilde\mu$ are equivalent, $\nu(x)=\nu(\bold1-x)$ for almost all $x\in X$  and there exist a measurable family of operators 
$$r(x)\:H_x \to H_{\bold 1 -x}\approx H_x$$
which are antilinear, preserve the norm and 
satisfy
$$
r(x)r(\bold1-x)=1, \qquad r(x)c_k(\bold 1-x)=(-1)^{k}c_k(x)r(x+\delta_k)\tag {3.6}$$
for all $k\in \Bbb N$ and almost all $x\in X$.
\endproclaim

\demo{Proof} 
Assume that  $S$ is an invariant real structure on $H$. To see that $\mu$ and $\tilde\mu$ are equivalent, it is enough to prove that $\tilde\mu(E)>0$ whenever  $E\subset X$  is measurable and $\mu(E)>0$. To this end, write $X=\cup_{n>1}E_n$ where $E_n=\{x\: \nu(x)=n\}$. We can assume that $E \subset E_n$  for some $n$. Let $h$ be a unit vector of $H_x$ for some $x\in E$. Define the functions $f(x)=\chi_E (x) h$ and $g(x)=(Sf)(x)$. As $S$ preserves norms 
$$
\mu(E)=\|f\|= \|Sf\| = \|g\|
$$
By the Lemma 3.4, the support of $g$ is $\bold 1 - E$, so
$$
\aligned
\mu(E)&=\left(\int_{\bold 1-E} (g(x),g(x))\,d\mu(x)  \right)^{\frac1 2} \\
	&=\left(\int_{E} (g(\bold 1 -x),g(\bold 1 -x))\,d\tilde\mu(x)  \right)^{\frac1 2} \\
\endaligned
$$
Then, $\tilde\mu(E)=0$ implies $\mu(E)=0$, which is a contradiction.

To see that $\nu(x)=\nu(\bold1-x)$ for almost all $x\in X$, assume the contrary. Then for some $n$ there exists a set $E\subset E_n$ such that $\mu(E)>0$ and $\bold 1 -E \subset E_{n-1}$. 
Let $\{h_i\}_{i=1}^n$ be an orthonormal basis of the $n$-dimensional Hilbert space $H_n$ ($H_x \approx H_n$  for all $x\in E$). Set $f_i(x)=\chi_E(x)h_i$ and $g_i=Sf_i$ for $i=1,\dots,n$. 
We will prove that for almost all $x$, $\left\{g_1(x),\dots, g_n(x)\right\}$  are $n$ non-zero orthogonal vectors in $H_{\bold1 -x}$, is a contradiction, since we have assumed $\bold 1 -E \subset E_{n-1}$.

For an arbitrary $F\subset E$ with $\mu(F)>0$, we have that
$$
\aligned
0\neq \|L_{\chi_{F}}f_i\|^2 &=\|SL_{\chi_{F}}f_i\|^2 =\|L_{\chi_{\bold1 -F}}Sf_i\|^2\\
				&=\int_{X}  \|\chi_{\bold1 -F}(x) g_i(x)\|^2\ d\mu(x)\\
				&=\int_{\bold1 -F} \|g_i(x)\|^2\ d\mu(x)
\endaligned
$$
Since $F$ is arbitrary and $\mu$ and $\tilde\mu$ are equivalent we have that $\|g_i(x)\|\neq 0$ a.e., and hence $g_i(x)\neq0$ a.e.

We now prove that $\nu(x)=\nu(\bold 1-x)$ almost everywhere. If not, there would be a set of positive measure $E$ such that $E\subset E_n$ and $\bold 1 -E \subset E_{n-1}$, for some $n$. Since $\mu$ and $\tilde\mu $ are equivalent, $\mu( \bold 1 -E)>0$. Since $\{f_i\}_{i=1}^n$ is an orthogonal set in $H$, so is $\{L_{\chi_F}f_i\}_{i=1}^n$ for any given subset $F\subset E$. By Lemma (3.1), the elements $SL_{\chi_F}f_i=L_{\chi_{\bold 1 -F}}g_i$ are also orthogonal. Therefore, for all  $i\neq j$
$$
\int_{\bold 1-F} (g_i(x),g_j(x)) d\mu(x) =0.
$$
Since $\mu\approx\tilde\mu$ and $F$ is arbitrary, we conclude that for almost all $x \in \bold 1-E$, 
$ (g_i(x),g_j(x)) =0.$ 
On the other hand, $\tilde\mu(F)= \int_{\bold 1-F} \|g_j(x)\|^2 d\mu(x)$   for almost all $x\in \bold 1-F$. Therefore $\|g_j\|\neq 0$, so that the $g_j$ cannot vanish. Since $E$ has dimension $n-1$, this is a contradiction. We conclude that $\nu(x)=\nu(1-x)$ for almost all $x$.

Up to identifications, we may now assume that 
$H_x=H_{\bold1-x}$. Define the operator $T: H\rightarrow H$ by
$$
(Tf)(x)= \sqrt{\frac{d\tilde\mu(x)}{d\mu(x)}} f(\bold 1-x)
$$
$T$ is $\Bbb C$-linear, unitary and satisfies the relations
$$
T^2=1  \qquad \qquad T N_k = N'_k T
$$
The product $ST$ is then an antilinear operator on $H$ which conmutes with all the $N_k$ and $N'_k$. Therefore it can be represented as multiplication by a measurable operator-valued function $r(x)$,
$$
(ST f)(x)= r(x)f(x)
$$
where $r(x): H_x\rightarrow H_x=H_{1-x}$ is antilinear and preserves norm. Since $T$ and $S$ are involutions, we obtain, respectively,
$$
(S f)(x)= r(x)(T f)(x) = r(x) \sqrt{\frac{d\tilde\mu(x)}{d\mu(x)}} f(\bold 1 -x).
$$
and
$$
r(x)r(\bold 1-x)=1.
$$
$S$ anticommutes with the operators  $A_k$ and $B_k$, because it commutes with $\Bbb J$. Therefore
$$
r(x)(-1)^{x_1+\dots+x_{k-1} +k-1} c_k(\bold 1-x)=(-1)^{x_1+\dots+x_{k-1}+1}c_k(x)r(x+\delta_k)
$$
and
$
r(x)c_k(\bold 1-x)=(-1)^{k}c_k(x)r(x+\delta_k)
$
follows.

Conversely, if $\mu$ and $\tilde\mu$ are equivalent, $\nu(x)=\nu(\bold1-x)$ a.e. and there exist antilinear operators $r(x)\:H_x \to H_x$ that preserve norm and satisfy (3.6) a.e., we can define $T$ and $S$ as above. Is straightforward to see that $S$ is in fact an invariant real form of $H$. 
\qed
\enddemo

\

If $\dim Z=2m <\infty$,  $\mu$ and $\tilde\mu$ are equivalent because they are quasi-invariant discrete measures concentrated in the same coset of $\Delta$ (namely, $X=\Delta$ itself) \cite {GW}. From (3.5) we deduce 
$$r(\bold1)=(-1)^{\frac {m(m+1)}2}r(\bold 0).$$ 
Assuming, as we may, that $r(\bold 0)$ is the standard conjugation on $\Bbb C$, we see that $H$ splits $/\Bbb R$ if and only if  $m(m+1)/2$ is an even integer, that is, for $m\equiv 0,3$  modulo $4$, as  is well known.

\

In the infinite-dimensional case there are always plenty of solutions $r(x)$ to the equations (3.6). They can all be described by fixing a full set $A\subset X$ of representatives of $X$ modulo $\Delta\cup{(\Delta+1)}$, defining $r(a)$ on $A$ arbitrarely and extending it to all of $X$ according to the rules 
$$r(\bold 1 - a)=r(a)^{-1}, \qquad r(x+\delta_k)=(-1)^{k}c_k(x)^{-1}r(x)c_k(\bold 1-x).$$
In this sense, any GW module with $\nu$ invariant and $\mu$  quasi-invariant under translation by $1$, ``splits over $\Bbb R$". But since the set $A$ is not measurable, the resulting operators $r(x)$ may not be either.  
The difficulty of deciding whether a given GW module actually splits, i.e., whether there is a measurable $r(x)$, depends on the module, as we see next. 
Assume now that $Z$ is infinite dimensional and separable.

\proclaim{Corollary 3.7} If  $\mu$ is discrete and $H$ is irreducible over $\Bbb C$, then it is irreducible over $\Bbb R$. In particular, this is the case for the Fermi-Fock representation.
\endproclaim

\demo{Proof} A proper invariant real subspace $U\subset H$ must be a real form of $H$. Indeed, $U\cap iU$ is complex, invariant and proper, hence 
equal to $\{0\}$ by irreducibility. By the same reason, $U\oplus iU$ is complex, invariant and not $\{0\}$, hence equal to $H$.
Since the measure $\mu$  is discrete and  $H$ is irreducible, $\mu$ is concentrated in some translate  $x_o+\Delta$, $x_o\in X$ \cite{GW}. As $X$ is infinite, $\bold 1-(x_o+\Delta)$ and $x_o+\Delta$ are disjoint. This implies that $\mu$ and $\tilde\mu$ are not equivalent and therefore   $H$ cannot have an invariant real form. 

The Fermi-Fock representation corresponds to the triple $(\mu,1,1)$, where $\mu(\delta_k)=1$ for all $k$ and $\mu(X-\Delta)=0$, hence it is supported in the discrete set $\Delta$.
\qed
\enddemo

\proclaim{Corollary 3.8} Let  $\mu_X$ be the Haar measure on $X$. Then the representation associated to the triple 
$(\mu_X,1,1)$ is irreducible over $\Bbb R$. 
\endproclaim

\demo{Proof} 
Irreducibility over $\Bbb C$ follows from the ergodicity of the Haar measure, as argued in [GW]. 
Irreducibility over $\Bbb R$ follows from Theorem 4 in [G], which implies that {\it any} function $r(x)$ satisfying
$$r(x+\delta_k)=(-1)^kr(x)$$
 must be non-measurable. 
\qed
\enddemo

\

While determining if a  GW module has an invariant real form may not always be easy, it is still possible to obtain all the real $C(H)$ modules up to orthogonal equivalence. The next and final Theorem shows that any pair consisting of a GW module together with an invariant real form, is unitarely equivalent to one of the following standard type.

Let 
$$K=\int_X^{\oplus}K_x \ d\mu(x)$$ be a direct integral of {\it real}  Hilbert spaces satisfying 
$$K_{x+\delta} = K_x,\qquad K_{1-x} = K_x$$
for all $\delta\in \Delta$ and almost all $x\in X$. With $H = {\Bbb C}\otimes K$ and $H_x={\Bbb C}\otimes K_x$,
$$H = \int_X^{\oplus}H_x \ d\mu(x).$$ 
Clearly,  $H_{x+\delta} = H_x=H_{1-x}$, $K_x$ is a canonical real form of $H_x$ and $K$ one of $H$. Denote by $\ \bar{}\ $ the corresponding conjugations. By definition, for $f\in H$,  
$$\bar f (x) = \overline{f(x)}.$$
This determines corresponding conjugations on  the endomorphisms $A$ of these spaces by $\bar A (f) = \overline{A(\bar f)}$.   

\proclaim{Theorem 3.9} If $(\mu,\nu,\{c_k\})$ is a GW module modeled on the Hilbert space $H$ with $$\tilde\mu\approx\mu,\qquad \nu(1-x)=\nu(x),\qquad \overline{c_k(1-x)} = (-1)^k c_k(x),$$
then
$$H^{\Bbb R} = \{f\in H: \overline{f(x)} = \sqrt{\frac{d\mu(1-x)}{d\mu(x)}} {f(1-x)}\}$$
is an invariant real form of $H$. Conversely, 
any pair consisting of a GW module together with an invariant real form, is unitarely equivalent to one of that type.
\endproclaim

\demo{Proof} 

For the first assertion we just need to apply Theorem (3.5) with $r(x)f(x) = \overline{f(x)}$.

For the converse, start with any split module and let $S,T,r(x),$ be as in the proof of (3.5). In it we have fixed  identifications $\psi_x: H_x\rightarrow H_{1-x}$ in order to define the operator $T$. Explicitely, letting $\Phi(x)=\sqrt{\frac{d\tilde\mu(x)}{d\mu(x)}}$, the ``actual" operator is $(\tau f)(x) = \Phi(x)f(1-x)$ and 
$$(Tf)(x)=  \psi_x^{-1}(\tau f)(x)=
\Phi(x)\psi_x^{-1}(f(1-x)).$$
If $\tilde \psi_x: H_x\rightarrow H_{1-x}$ is another identification, the corresponding operator will be
$$ (\tilde T f)(x) = \tilde \psi_x^{-1}((\tau f)(x))
= \tilde \psi_x^{-1}\circ \psi_x \circ \psi_x^{-1}((\tau f)(x)) = h_x ((Tf)(x))$$
with $h_x = \tilde \psi_x^{-1}\circ \psi_x$ unitary on $H_x$.  

Define $u(x)=H_x\rightarrow H_x$ by
$$u(x)f(x) = r(x)\overline{ f(x)}.$$
These are $\Bbb C$-linear and norm-preserving, hence unitary and measurable in $x$, because $r(x)$ is so. Change the chosen identification $H_x\approx H_{1-x}$ by $u(x)$, so that the new operator $T$ will be 
$$\tilde Tf(x)= u(x)^{-1}(Tf)(x).$$
For this identification and the equivalent real structure
$$\tilde S = u^{-1}Su,$$
the function $\tilde r(x)$ is
$$ \tilde r(x)f(x) = \tilde S \tilde T f(x)= u^{-1}Su 
u(x)^{-1}(Tf)(x) = u(x)^{-1}STf(x) = u(x)^{-1}r(x)f(x) = 
\overline{f(x)}.$$
This effectively shows that, up to equivalence, we may take $r(x)$ to be conjugation with respect to any basic flat field of real forms. 

Since the operators $c_k(x)$ do not depend on the identification $H_x=H_{1-x}$, the relations (3.6) must still hold for the new $\tilde r(x)$, which translates into $\overline{c_k(1-x)} = (-1)^k c_k(x)$.

\qed
\enddemo

\

As an application, let us describe the simplest irreducible representation of $C_\Bbb C(Z)$ that splits over $\Bbb R$., i.e., as a representation of $C(Z)$. The two irreducible summands will be equivalent and truly real, in the sense that they do not arise from any complex representation of $C_\Bbb C(Z)$ by restriction of the scalars.

Let $\mu=\mu_X$ be the Haar measure on $X$ and $\nu=1$, so that $H_x=\Bbb C$, $K_x=\Bbb R$, 
$$H=L_2(X)$$ 
(dropping the reference to the Haar measure) and the conjugation $\bar{ } $ is given by $\bar f (x)=\overline{f(x)}$. As $c_k$'s, which are now functions from $X$ to the circle $\Bbb T$, choose 
$$\eqalign{c_{2\ell}(x)&=1\cr
c_{4 \ell+1}(x) &= (-1)^{x_{4\ell+3}}\cr
c_{4\ell+3}(x) &= (-1)^{x_{4\ell+1}}\cr}.$$

\proclaim{Theorem 3.11} These functions satisfy (2.3) and, therefore, (2.5) defines a representation of $C_\Bbb C(Z)$ in $L_2(X)$. This representation is irreducible. The  real subspace
$$
L_2(X)^{\Bbb R} = \{f\in L_2(X): \ f(1-x)= \overline{f(x)} \}
$$
is invariant and irreducible under $C(Z)$. Furthermore, the  representation so obtained does not arise from any representation of $C_\Bbb C(Z)$ by restriction of the scalars. 
\endproclaim

\demo{Proof} The irreducibility of the corresponding GW module $H$ follows from the ergodicity of the Haar measure, as in (3.8).  

It is straightforward to check that the functions $c_k$ satisfy (2.3) and (3.10), so, by Theorem (3.9), the corresponding $\Bbb J$ must leave the real form $H^{\Bbb R}$ invariant. Of course, this can be deduced by direct calculation as well, using (2.5). 
$H^{\Bbb R}$ must be irreducible under $C(Z)$, since any closed invariant subspace generates a closed $C_{\Bbb C}(Z)$-invariant subspace in $H$. 

Finally, suppose that the representation of $C(Z)$ in $H^{\Bbb R}$ could be extended to one of $C_{\Bbb C}(Z)$ (in $H^{\Bbb R}$ itself). Denote by $T$ the operation representing multiplication by $\sqrt{-1}$: $T$ is an
orthogonal complex structure in $H^{\Bbb R}$ commuting with $C(Z)$. Its unique $\Bbb C$-linear extension to all of $H=H^{\Bbb R}\oplus iH^{\Bbb R}$ is unitary and  commutes with all the $\Bbb J$. As we have already mentioned, this implies that $T$ is given pointwise, by an operator-valued measurable function $k(x)$:
 $(Tf)(x) = k(x)f(x).$ In the present case, $k(x)$ is complex valued. Since $k(x)^2=-1$, we can write it as 
$$k(x)=\epsilon(x)i$$
for some measurable $\epsilon: X\rightarrow \{\pm1\}$. The condition for  $T$ to leave invariant the real form $H^{\Bbb R}$ and to commute with the representation $\Bbb J$ amount to, respectively, 
$$
\epsilon(1-x)=-\epsilon(x)\qquad\qquad \epsilon(x+\delta_k) = \epsilon(x)
$$
for almost all $x$ and all $k$. The second equation implies that $\epsilon$ is actually constant on each $\Delta$-equivalence class. [According to [Go] ???], $\epsilon$ must then be constant almost everywhere, contradicting the first equation. Hence, no such $T$ can exist.
\qed
\enddemo

Representations with these parameters 
can be realized on the standard $L_2(\Bbb T)$ of complex-valued functions on the circle since, as a measure space,  $(X,\mu_X)$ is the same as the interval $(0,1)$ -hence to the circle, equipped with the Lebesgue measure and $\dim H_x=1$. In this identification, translations in  $X$ do not correspond to rigid rotations of $\Bbb T$. However, the operation $x\mapsto 1-x$ in $X$ (changing all the elements of $x$) does correspond to $x\mapsto 1-x$ on $(0,1)$ which, on $\Bbb T\subset \Bbb C$, becomes ordinary complex conjugation. 
The operators $c_k$ can then be viewed as functions
$$\tilde c_k: {\Bbb T}\rightarrow \Bbb T,$$
satifying, of course, (2.3).
According to Theorem 3.9, those that split over $\Bbb R$ can be realized with 
$$
L_2(\Bbb T)^{\Bbb R} = \{f\in L_2(\Bbb T): \overline{f(t)} = {f(\bar t)}\}.
$$
as the invariant real form. The condition on the functions $c_k$ for this real form to be invariant under the $\Bbb J$'s  becomes 
$$\overline{c_k(t)} = (-1)^k c_k(\bar t).\leqno{(3.10)}$$
for $t\in\Bbb T$.

\enddocument
\end